\documentclass[11pt]{article}
\input amssym.def
\input amssym.tex
\usepackage{epsfig,psfig}
\def\binom#1#2{{#1}\choose{#2}}

\def\slfrac#1#2{\hbox{\kern.1em %
 \raise.5ex\hbox{\the\scriptfont0 #1}\kern-.11em %
 /\kern-.15em\lower.25ex\hbox{\the\scriptfont0 #2}}}

\newcommand{\eqn}[1]{(\ref{#1})}
\newcommand{\hsp}{\hspace*{\parindent}}

\newcommand{\eeq}{\end{equation}}
\newcommand{\beql}[1]{\begin{equation}\label{#1}}

\newcommand{\FF}{{\Bbb F}}

\newcommand{\sC}{{\cal C}}

\makeatletter
\def\@sect#1#2#3#4#5#6[#7]#8{\ifnum #2>\c@secnumdepth
     \def\@svsec{}\else
     \refstepcounter{#1}\edef\@svsec{\csname the#1\endcsname.\hskip .75em }\fi
     \@tempskipa #5\relax
      \ifdim \@tempskipa>\z@
        \begingroup #6\relax
          \@hangfrom{\hskip #3\relax\@svsec}{\interlinepenalty \@M #8\par}%
        \endgroup
       \csname #1mark\endcsname{#7}\addcontentsline
         {toc}{#1}{\ifnum #2>\c@secnumdepth \else
                      \protect\numberline{\csname the#1\endcsname}\fi
                    #7}\else
        \def\@svsechd{#6\hskip #3\@svsec #8\csname #1mark\endcsname
                      {#7}\addcontentsline
                           {toc}{#1}{\ifnum #2>\c@secnumdepth \else
                             \protect\numberline{\csname the#1\endcsname}\fi
                       #7}}\fi
     \@xsect{#5}}
\def\@begintheorem#1#2{\it \trivlist \item[\hskip \labelsep{\bf #1\ #2.}]}

\def\plain{plain}\ifx\fmtname\plain\csname fi\endcsname
     
     \input docstrip
     \preamble

     Do not distribute the stripped version of this file.
     The checksum in the header refers to the documented version.

     \endpreamble
     \generateFile{here.sty}{t}{\from{here.doc}{}}
     \endinput
\fi
\ifcat a\noexpand @\let\next\relax\else\def\next{%
    \documentstyle[here,doc]{article}\MakePercentIgnore}\fi\next
\ifx\@Hxfloat\@Hundef\else\expandafter\endinput\fi
\let\@Hxfloat\@xfloat
\def\@xfloat#1[{\@ifnextchar{H}{\@HHfloat{#1}[}{\@Hxfloat{#1}[}}
\def\@HHfloat#1[H]{%
\expandafter\let\csname end#1\endcsname\end@Hfloat
\vskip\intextsep\vbox\bgroup\def\@captype{#1}\parindent\z@
\ignorespaces}
\def\end@Hfloat{\egroup\vskip \intextsep}

\makeatother

\begin{document}
\begin{center}
{\Large {\bf On Kissing Numbers in Dimensions 32 to 128}} \\
\vspace{1.5\baselineskip}
{\em Yves Edel} \\
Mathematisches Institut der Universit\"{a}t \\
Im Neuenheimer Feld 288 \\
69120 Heidelberg, Germany \\
\vspace{1\baselineskip}
{\em E. M. Rains} and {\em N. J. A. Sloane} \\
\vspace*{.5\baselineskip}
AT\&T Labs-Research \\
180 Park Avenue \\
Florham Park, NJ 07932-0971, USA \\
\vspace{1.5\baselineskip}
July 19 1998 \\
\vspace{1.5\baselineskip}
{\bf ABSTRACT}
\vspace{.5\baselineskip}
\end{center}
\setlength{\baselineskip}{1.5\baselineskip}

An elementary construction using binary codes gives new record
kissing numbers in dimensions from 32 to 128.
\clearpage
\thispagestyle{empty}
\setcounter{page}{1}

\section{Introduction}
\hsp
Let $\tau_n$ denote the maximal kissing number in dimension $n$, that is, the greatest number of $n$-dimensional spheres that can touch another sphere of the same size.
Although asymptotic bounds on $\tau_n$ are known \cite{SPLAG}, little is known about explicit constructions, especially for $n > 32$.
Up to now the best explicit constructions have come
from lattice packings.
The kissing number $\tau$ of the
Barnes-Wall lattice\footnote{The subscript gives the dimension.}
$BW_n$ in dimension $n=2^m$ is $\prod_{i=1}^m (2^i +2)$,
although for $m \ge 5$ this is weak (146,880, 9,694,080 and 1,260,230,400 in dimensions 32, 64 and 128,
for example).
In contrast, Quebbemann's lattice
$Q_{32}$ \cite{Que6}, \cite[Chap.~8]{SPLAG} has $\tau ={} $261,120.

In recent years the kissing numbers of a few other lattices in dimensions $> 32$ have been determined.
Nebe \cite{Nebe98} shows that the Mordell-Weil lattice $MW_{44}$ has
$\tau ={} $2,708,112.
Nebe \cite{Nebe98b} shows that a 64-dimensional lattice constructed
in \cite{Nebe98} is extremal 3-modular, and so by modular form theory has
$\tau ={} $138,458,880.
Bachoc and Nebe \cite{BacoN96a} give an 80-dimensional lattice with $\tau ={} $1,250,172,000.
Elkies \cite{Elki98} calculated the
kissing number of his lattice $MW_{128}$:
it is 218,044,170,240, over 170 times that of $BW_{128}$.

In the present note, we show that an elementary construction using binary codes gives better values than all of these.
However, our packings are just local arrangements of spheres around the origin:
we do not know if they can be modified to produce dense {\em infinite} packings.

\section{The construction}
\hsp
Let $\sC(n, d)$ (resp. $\sC (n, d, w)$) denote a set of binary vectors of length $n$ and Hamming distance
$\ge d$ apart (resp. and with constant weight $w$).
The maximal size of such a set is denoted by $A(n, d)$ (resp. $A(n, d, w)$)
\cite{BrSSS}, \cite{MS}.

One way to achieve the kissing number $\tau_8 = 240$ in eight dimensions is to take as centers of spheres the vectors of shape $\pm 1^8$, with a unique
support (a code $\sC(8, 8, 8)!$) and signs taken from a $\sC (8, 2)$, together
with the vectors of shape $\pm 2^2 0^6$, where the supports are taken from a
$\sC(8, 2, 2)$ and the signs from a $\sC (2,1)$.
Taking all these codes to be as large as possible, we obtain a total of
$$A(8, 8, 8) A(8, 2) + A(8, 2, 2) A(2, 1) =
1 \cdot 2^7 + {\binom{8}{2}}2^2 = 240
$$
spheres touching the sphere at the origin.

Our construction generalizes this as follows.
For a given dimension $n$, we choose a sequence of support sizes
$n_0$, $n_1, \ldots, n_\mu$ satisfying
\beql{Eq1}
n \ge n_0 \ge 4n_1 \ge 4^2 n_2 \ge \cdots \ge
4^\mu n_\mu \ge 1 ~.
\eeq
The $\nu^{\rm th}$ set of centers that we use, for $0 \le \nu \le \mu$,
consists of vectors of shape $\pm a_\nu^{n_\nu} 0^{n-n_\nu}$, where
$a_\nu = \sqrt{n_0 / n_\nu}$, the supports are taken from a
$\sC (n, n_\nu , n_\nu )$ and the signs from a $\sC(n_\nu,  \left\lceil \frac{n_\nu}{4} \right\rceil)$.
With optimal codes, the total number of vectors is
\beql{Eq2}
\sum_{\nu =0}^\mu A(n, n_\nu , n_\nu ) A\left( n_\nu , 
\left\lceil \frac{n_\nu}{4} \right\rceil \right) ~.
\eeq
It is easy to check that all vectors have length $\sqrt{n_0}$, and that by \eqn{Eq1} the distance between any two distinct vectors is $\ge \sqrt{n_0}$.
It follows that \eqn{Eq2} is a lower bound on $\tau_n$.
\subsection*{Remarks}
\hsp
(1) Even if we do not know the exact values of $A(n, d, w)$ and $A(n, d)$ mentioned in \eqn{Eq2}, we can replace them by any available lower bounds, and still obtain a lower bound on the kissing number $\tau_n$.
There is some freedom in choosing the $n_\nu$, which helps to compensate for our ignorance.

(2) A table of lower bounds on $A(n, d)$ has been given by Litsyn \cite{Lits98},
extending the table in \cite{MS}.
A table of lower bounds on $A(n, d, w)$ for $n \le 28$ is given in \cite{BrSSS},
but for larger $n$ little is known.
A very incomplete table for $n > 28$ can be found in \cite{RaSl98b}.

(3) The construction gives a set of points on a sphere with angular separation of $60^\circ$.
It can obviously be modified to produce spherical codes with other angles.

\section{Examples}
\hsp
We illustrate the construction by giving new records in dimensions 32, 36, 40, 44, 64, 80 and 128.
For other examples see \cite{NeSl}, and for further details about the codes see
\cite{Lits98}, \cite{RaSl98b}.

\paragraph{n=32.}
We take $n_0 = 32$, $n_1= 8$, $n_2 =2$ and use 
$A(32, 8) \ge 2^{17}$ from \cite{Che3},
$A(32, 8, 8) \ge 1117$ from the complement of a lexicographic code $\sC (32, 8, 24)$
(cf. \cite{Con41}), obtaining a kissing number of $A(32, 32, 32) A(32, 8) + A(32, 8, 8) A(8, 2) + A(32, 2, 2) A(2, 1) \ge 1 \cdot 2^{17} + 1117 \cdot 2^7 + {\binom{32}{2}}  \cdot  2^2 ={} $276,032.

\paragraph{n=36.}
Let the 36 coordinates be labeled $(i, j)$, $0 \le i,j \le 5$, and let the symmetric group $S_6$ act by $(i, j) \to (i^\pi , j^{\sigma (\pi )} )$, where $\pi \in S_6$ and $\sigma$ is the outer automorphism of $S_6$.
One can find a set of 17 orbits under the alternating group $A_6$,
of sizes ranging from 45 to 360, whose union forms a constant weight code
showing that
$A(36, 8, 8) \ge 2385$.
We take $n_0 = 32$, $n_1 =8$, $n_2 =2$ and obtain a kissing number of
$A(36, 32, 32) A(32, 8) + A(36, 8, 8)A(8, 2) + A(36, 2, 2) A(2, 1) \ge 1 \cdot 2^{17} + 2385 \cdot 2^7 + {\binom{36}{2}}  \cdot 4 ={} $438,872.

An alternative approach can be based on Warren D. Smith's discovery (personal communication, May 1997) that the 2754 minimal vectors of the self-dual
length 18 distance 8 code over $\FF_4$ \cite{Mac4} yields
$\tau_{36} \ge 2754 \cdot 2^7 ={} $352,512 by changing any even number of signs.
By adjoining additional vectors with fractional coordinates R.~H. Hardin and N.~J.~A. Sloane increased this to 386,570, which held the record until it
was overtaken by the present construction.
It is quite possible that with better clique-finding the $\FF_4$ approach will regain the lead.

\paragraph{n=40.}
We take $n_0 = 40$, $n_1 =8$, $n_2=2$, use a lexicographic code for $A(40, 8, 8)$, and obtain $A(40, 40, 40) A(40, 10) + A(40, 8, 8) A(8, 2) + A(40, 2, 2) A(2, 1) \ge 1 \cdot 589824 + 3116 \cdot 2^7 + {\binom{40}{2}}  \cdot  2^2 ={} $991,792.

\paragraph{n=44.}
$A(44, 44, 44) A(44, 11) + A(44, 8, 8)A(8, 2) + A(44, 2, 2) A(2, 1) \ge 1 \cdot 2^{21} + 6622 \cdot 2^7 + {\binom{44}{2}}  \cdot 4 ={} $2,948,552.

\paragraph{n=48.}
In 48 dimensions the three known extremal unimodular lattices
\cite{SPLAG}, \cite{Nebe98} have kissing number 52,416,000.
Our present construction gives less than half this value.

\paragraph{n=64.}
The words of weight 16 in an extended cyclic code $\sC(64, 16)$ of size $2^{28}$ from \cite{Pet3} show that $A(64, 16, 16) \ge {} $30,828.
In this way we obtain a kissing number of 331,737,984.

\paragraph{n=80.}
By taking 4 orbits under $L_2(79)$ we obtain
$A(80, 16, 16) \ge {} $143,780.
We take $n_0 = 64$, $n_1=16$, $n_2 =4$, $n_3 =1$ and
obtain $\tau \ge  {}$1,368,532,064.

\paragraph{n=128.}
This is the most dramatic improvement, so we give a little more detail.
Our construction uses:
$$
\begin{array}{lll}
A(128, 128, 128) A(128, 32) ~{\rm vectors}~ \pm 1^{128}: & \ge & 1 \cdot 2^{43} \\
A(128, 32, 32)  A(32, 8) ~{\rm vectors}~ \pm 2^{32} 0^{96}: & \ge & 512064  \cdot 2^{17} \\
A(128, 8, 8) A(8, 2) ~{\rm vectors}~ \pm 4^{8} 0^{120}: & \ge & 2704592 \cdot 2^7 \\
A(128, 2, 2) A(2, 1) ~{\rm vectors}~ \pm 8^{2} 0^{126}: & \ge & {\binom{128}{2}}  \cdot 4 \\ \cline{3-3}
\multicolumn{2}{r}{\mbox{for a total of}} & \mbox{8,863,556,495,104}
\end{array}
$$
Here $A(128, 32) \ge {} 2^{43}$ comes from a BCH code \cite[p.~267]{MS},
$A(128, 32, 32) \ge {} $512064 from a union of two orbits under $L_2(127)$,
$A(32, 8) \ge {} 2^{17}$ from \cite{Che3},
and $A(128, 8, 8) \ge {} $2704592
is obtained by shortening a  $\sC (129, 8, 8)$ of size 2883408
formed from the union of 11 orbits of size 262128 under $L_2(128)$.
The result is more than 40 times that of the Mordell-Weil lattice.

\vspace*{+.3in}

We do not expect any of these new records to survive for long, since our lower bounds for $A(n, d)$ and $A(n, d, w)$ are very weak.
However, it is interesting that such a simple construction gives such dramatic improvements over the kissing numbers of the best lattices known.

\vspace*{+.3in}
{\bf Postscript.} Victor Zinoviev has pointed out to us that 
in 1992 he and T. Ericson \cite{EZ92} obtained kissing numbers of 
858800 in dimension 40 and 273935235 in dimension 64
by similar methods.
Warren D. Smith (unpublished) used these methods to
construct spherical codes in the 1980's.

\end{document}